\documentclass[11pt]{amsart}
\usepackage{amsmath}
\numberwithin{equation}{section}
\newtheorem{thm}{Theorem}[section]
\newtheorem{prop}[thm]{Proposition}
\newtheorem{lemma}[thm]{Lemma}
\newtheorem{cor}[thm]{Corollary}
\theoremstyle{remark}
\newtheorem*{remark}{Remarks}
\newcommand{\N}{{\mathbb N}}
\newcommand{\R}{{\mathbb R}}

\newcommand{\Z}{{\mathbb Z}}

\title[Spectrum and Embeddings]{The Spectrum and Isometric Embeddings 
of Surfaces of Revolution \\}
\author{Martin Engman}
\begin{document}
\address{Departamento de Ciencias y Tecnolog\'{i}a, Universidad Metropolitana, 
San Juan, PR 00928}
\email{mengman@excite.com}
\thanks{Partially supported by the NSF Grant: Model Institutes for Excellence at 
UMET and PR-EPSCoR}
\subjclass{Primary 58G25, 58G30, 53C20;
Secondary  35P15}
\dedicatory{For Gus and Sonia}
\begin{abstract}
{An upper bound on the first $S^{1}$ invariant eigenvalue of the Laplacian 
for $S^1$ invariant metrics on $S^2$ is used to find obstructions to the 
existence of isometric embeddings of such metrics in $(\R^3, can)$. As a 
corollary we prove: If the first four distinct eigenvalues have even 
multiplicities then the surface of revolution cannot be isometrically 
embedded in $(\R^3, can)$. This leads to a generalization of
a classical result in the theory of surfaces.}  
\end{abstract}
\maketitle
\section{Introduction}

The problem of isometrically embedding $(S^2,g)$ in $(\R^3, can)$ has a 
long history which goes back at least as far as 1916. In that 
year, Weyl, \cite{we}, 
and in the years since, Nirenberg, \cite{ni}, Heinz, \cite{he}, Alexandrov, 
\cite{al}, and Pogorelov, \cite{po}, to name a few, proved embedding 
theorems of various orders of differentiability in case the Gauss curvature 
is positive. 
A recent result of Hong and Zuily \cite{hz} addresses the case of 
non-negative curvature. But, of course, not every metric on $S^2$ admits such
an isometric embedding. The reader may refer to Greene, \cite{gr}, 
wherein one finds examples of smooth metrics on 
$S^2$ for which there is no $C^2$ isometric embedding in $(\R^3, can)$.

In the presence of examples such as Greene's, one might naturally ask if 
there exist intrinsic geometric conditions on metrics which obstruct such 
isometric embeddings.
Inasmuch as the above mentioned embedding theorems require, at least, 
non-negativity of the Gauss curvature, one must look for embedding obstructed
metrics among those with some negative curvature. Of course, having some
negative curvature is not enough, but one might hope that some stronger
condition, associated with the existence of negative curvature at a point, 
might satisfy our requirements. 
The purpose of this paper is to provide, in a special case, 
conditions on the spectrum of the Riemannian manifold which are intrinsic 
obstructions to the above isometric embedding problem. 

It is no surprise that the spectrum might make an appearance in this subject.
There is an extensive literature which associates the spectrum with (more 
generally) isometric immersions (see, for example, \cite{chen} or \cite{li} 
among others). Much of this work relates the spectrum to the mean curvature,
and is associated with the Willmore conjecture. By way of comparison, the
embedding problem of this paper is almost trivial, but it is the new and more 
intrinsic relation between the spectrum and embeddability which we hope the
reader will find interesting. 

In the case of $S^1$ invariant metrics on $S^2$ 
(i.e. surfaces of revolution), one can prove that, while the first eigenvalue
must be bounded above by $8\pi/area$ (Hersch's theorem, \cite{her}), 
the first $S^1$ invariant
eigenvalue can be arbitrarily large. At the same time, however, there is an
upper bound, depending on the metric, for the first $S^1$ invariant 
eigenvalue. We will prove that it is this upper bound which, upon exceeding
a certain critical value, becomes an obstruction to isometric embedding into
$(\R^3, can)$ (Necessarily, the same condition ensures that there is some 
negative curvature). As a result, if the first $S^1$ invariant eigenvalue 
becomes 
too large then the surface cannot be isometrically embedded in $(\R^3, can)$.
(See also Abreu and Freitas \cite{abr}.)

Another characteristic of the spectrum of a surface of revolution is that the
eigenspaces are even dimensional unless the eigenvalue happens to correspond
to an 
$S^1$ invariant eigenvalue. As a result, one way of increasing the first
 $S^1$ invariant eigenvalue is to insist that the multiplicities be even up
to a certain point. This leads to a result, proved in Section 5, that even 
multiplicity for the first 4 distinct eigenvalues is an obstruction to 
isometric embeddability.

In the last section, we will remark on how these results give a 
generalization of a well known corollary of the Gauss-Bonnet theorem 
regarding the existence of points with negative curvature.
 
The author is indebted to Andrew Hwang for a brief but enlightening 
conversation about momentum coordinates.  

\section{Isometric embeddings and Momentum Coordinates}

First, we will discuss a formulation of the condition on the metric which 
ensures that an $S^1$ invariant metric on $S^2$ may be isometrically
 embedded into $(\R^3, can)$. This condition is well known and quite 
elementary. The reader will find our treatment to be essentially equivalent 
to that of Besse \cite{be}, p. 95-105. 

Let $(M,g)$ be an $S^1$ invariant Riemannian manifold which is diffeomorphic
to $S^2$. We will assume the metric to be
$C^{\infty}$.
 Since $(M,g)$ has an effective $S^{1}$ isometry group there are exactly two 
fixed points. We call the fixed points $np$ and $sp$ and let $U$ be the
 chart $M \setminus \{np,sp\}$. On $U$ the metric has the form 
$ds \otimes ds + a^{2}(s)d\theta \otimes d\theta$ where $s$ is the 
arclength along a geodesic connecting $np$ to $sp$ and $a(s)$ is a 
function $a:[0,L]\rightarrow \R^+$ satisfying $a(0)=a(L)=0$ and 
$a'(0)=1=-a'(L)$.

It is easy to see that isometric embeddings of such metrics can be 
parametrized as follows:

\begin{equation}
\left\{ \begin{array}{ccl}
  \psi^1 & = & a(s)\cos \theta \\
  \psi^2 & = & a(s)\sin \theta \\
  \psi^3 & = & \pm \int^s \sqrt{1-(a')^2(t)}dt
                      \end{array} \right. \label{eq:param}
\end{equation} 
It is evident from this formula that $(M,g)$ can be isometrically $C^1$
 embedded
 in $(\R^3, can)$ if and only if

\begin{equation}
|a'(s)|\leq 1 \hspace{.1in} \mbox{for all} \hspace{.1in} s \in [0,L] . 
\label{eq:leq1}
\end{equation}

We will find it convenient to make a change of variables to so called 
\underline{momentum coordinates} (See Hwang and Singer \cite{hw}).
These are given by a diffeomorphism 
$(s,\theta) \rightarrow (x,\theta)$ where 
$x\equiv \phi :[0,L] \rightarrow [-1,1]$ is defined by:

\begin{equation}
x\equiv \phi(s)\equiv \int^s_c a(t)dt. \label{eq:phi}
\end{equation}
   
If we let $f(x) \equiv (a^2\circ \phi^{-1})(x)$, then in the new 
coordinates the metric on the chart $U$ takes the form

\begin{equation}
 g=\frac{1}{f(x)}dx \otimes dx + f(x)d\theta \otimes d\theta    
\label{eq:feq}
\end{equation}
where $(x,\theta) \in(-1,1) \times [0,2\pi)$. In these coordinates the
 conditions at the endpoints translate to $f(-1)=0=f(1)$ and 
$f'(-1)=2=-f'(1)$. 
In this form, it is easy to see that this metric has area $4\pi$ and that its 
Gauss curvature is given by 
$K(x) =(-1/2)f^{''}(x)$. It is also worth observing that the function $f(x)$
is the square of the length of the Killing field (infinitesimal isometry) 
$\partial/\partial \theta$ on the chart $U$.
 The canonical (i.e. constant curvature) metric is obtained by taking 
$f(x)=1-x^{2}$. 

Now using \eqref{eq:leq1} and the definition for $f$ it is a simple 
exercise from calculus to prove:

\begin{prop} Let $(M,g)$, with metric $g$ as in \eqref{eq:feq}, be 
diffeomorphic to $S^2$. $(M,g)$ can be isometrically $C^1$ embedded 
in $(\R^3,can)$ if and only if $|f'(x)| \leq 2$ for all $x \in [-1,1]$. 
\qed
\end{prop}

We will end this section with the comment that it is also very easy to see 
that (in our special case) non-negative curvature implies isometric 
embeddability since $K(x) \geq 0$ implies that $f'(x)$ is a non-increasing
function on $[-1,1]$ with maximum $2$ and minimum $-2$. One can find 
essentially the same comment on p-106 of \cite{be}.

\section{Some properties of the spectrum}

In the interest of presenting a self-contained exposition we will review 
some of the relevent 
facts about the spectrum (eigenvalues) of a surface of revolution in this
section. The interested reader may consult \cite{eng1}, \cite{eng2} and 
\cite{eng3} for further details.

Let $\Delta$ denote 
the scalar Laplacian on a surface of revolution $(M,g)$, where $g$ is 
given by \eqref{eq:feq} and let $\lambda$ be any eigenvalue of $-\Delta$.
We will use the symbols $E_{\lambda}$ and $\dim E_{\lambda}$ to denote the 
eigenspace for $\lambda $ and it's multiplicity respectively. In this paper 
the symbol $\lambda_m$ will always mean the $m$-th \underline{distinct} 
eigenvalue. We adopt the convention $\lambda_0=0$. 
Since $S^1$ (parametrized here by $0 \leq \theta <2\pi $) acts
 on
$(M,g)$ by isometries and because $\dim E_{\lambda_m} \leq 2m+1$ (see 
\cite{eng3} for the proof), the 
orthogonal decomposition of $E_{\lambda_m}$ has the special form

\begin{equation*}
E_{\lambda_m}= \bigoplus_{k=-m}^{k=m} e^{ik\theta}W_k
\end{equation*}
in which $W_k (=W_{-k})$ is the ``eigenspace" (it might contain only $0$)
 of the ordinary differential 
operator

\begin{equation*}
L_{k}=-\frac{d}{dx}\left(f(x)\frac{d}{dx}\right) + \frac{k^2}{f(x)}
\end{equation*}
with suitable boundary conditions. It should be observed that $\dim W_k \leq
1$, a value of zero for this dimension occuring when $\lambda_m \not\in Spec
L_k$.

It is easy to see that $Spec(-\Delta)= \bigcup_{k\in \Z} Spec L_k$, 
consequently the
spectrum of $-\Delta$ can be studied via the spectra 
$Spec L_{k}=\{0<\lambda_{k}^{1} < \lambda_{k}^{2} < \cdots 
<\lambda_{k}^{j} < \cdots \} \forall k \in \Z$. 
The eigenvalues $\lambda^j_0$ in the case $k=0$ above are called the 
\underline{$S^1$ invariant eigenvalues} since their eigenfunctions are 
invariant under the the $S^1$ isometry group. If $k \neq 0$ the eigenvalues
are called \underline{$k$ equivariant} or simply \underline{of type $k \neq 0
$}. Each $L_{k}$ has a Green's operator, 
$\Gamma_{k}:(H^{0}(M))^{\perp} \rightarrow L^{2}(M)$, whose spectrum 
is $\{ 1/\lambda_{k}^{j} \}_{j=1}^{\infty}$, and whose trace is defined by, 
$tr \Gamma_{k} \equiv \sum 1/\lambda_{k}^{j}$.

\begin{prop}[See \cite{eng1} and \cite{eng2}] With the notations as above:  

\begin{description}
   \item[i] $$tr \Gamma_{k}= \left\{ \begin{array}{cl}
   \frac{1}{2} \int_{-1}^{1} \frac{1-x^{2}}{f(x)} dx  & 
 \mbox{if $k=0$} \\
   \frac{1}{|k|} &  \mbox{if $k \neq 0$}
                     \end{array} \right. .$$
\item[ii] If $area(M,g)=4 \pi$ and
${\displaystyle \sum^{\infty}_{j=1} \frac{1}{\lambda^{j}_{0}} \leq
 \frac{\pi^2}{16}}$
then there exist points $p \in M$ such that $K(p)<0$.
  
\item[iii] For all $k \in \Z$ and $j \in \N $, 
$ \lambda_{k}^{j}=\lambda_{-k}^{j}$.
   \item[iv] $ \forall k \geq 1$ and $ \forall j \geq 0$, 
 $\lambda_{k+j} \leq \lambda_{k}^{j+1}$; and  $\lambda_1 \leq \lambda_0^1$.
 \item[v] $\dim E_{\lambda_{m}}$ is odd if and only if $\lambda_{m}$ is 
an $S^{1}$ invariant eigenvalue. \qed
 \end{description}
\end{prop} 

\begin{remark}

1.) One must be careful with the definition of $tr \Gamma_0$ since 
$\lambda_0 = 0 \in 
Spec L_0$. To avoid this difficulty we studied the $S^1$ invariant spectrum
of the Laplacian on 1-forms in \cite{eng2} and then observed that the 
non-zero eigenvalues are the same for functions and 1-forms.

2.) A slight modification of the proof of Proposition 3.1 ii.) 
(in \cite{eng2}) reveals that ${\displaystyle \sum^{\infty}_{j=1} 
\frac{1}{\lambda^{j}_{0}}
 \leq
 \frac{\pi^2}{16}}$ implies that $(M,g)$ cannot be isometrically embedded in
 $(\R^3, can)$. The reader will find that the results of this paper are an
improvement on this idea.
\end{remark}

\section{A sharp upper bound for the first eigenvalue}

In \cite{eng3} we derived sharp upper bounds for all of the distinct 
eigenvalues on a surface of revolution diffeomorphic to $S^2$. These 
estimates were obtained using the the $k$-type eigenvalues for $k \neq 0$.
In this section we will obtain a sharp bound for $\lambda_1$ using the 
$S^1$ invariant spectrum. In contrast with the more general result of Hersch 
\cite{her}, the reader will find that this bound exhibits, more explicitly, 
its 
dependence on the metric. This fact will play an important r\^{o}le in 
embedding problems.   

\begin{prop} 
Let $(M,g)$ be an $S^1$ invariant Riemannian manifold of area $4\pi$ which is
diffeomorphic
to $S^2$, with metric \eqref{eq:feq}. Let $ \lambda^1_0 $ be the first 
non-zero $S^1$
invariant eigenvalue for this metric, then

\begin{equation*}
\lambda^1_0 \leq \frac{3}{2} \int^{1}_{-1}f(x)dx
\end{equation*}

and equality holds if and only if $(M,g)$ is isometric to $(S^2,can)$.

\end{prop}

\begin{proof} The minimum principle associated with the first
 $S^1$ invariant eigenvalue problem,

\begin{equation}
L_0 u= -\frac{d}{dx}\left(f(x)\frac{du}{dx}\right)=\lambda^1_0 u, 
\label{eq:de}
\end{equation}
states that
\begin{equation}
 \lambda^1_0 \leq \frac{\int^{1}_{-1}f(x)(\frac{du}{dx})^2 dx}
{\int^{1}_{-1}u^2 dx} \label{eq:min}
\end{equation}
for all $S^1$ invariant functions $u \in C^{\infty}(M)$ with 
$u\perp \ker L_0$. Equality holds if and only if $u$ is an eigenfunction 
for $\lambda^1_0$. 
Since
$ \ker L_0$ consists of constant functions and  
$\int^{1}_{-1}x \cdot 1 dx=0$, we see that $u(x)=x$ is an admissible 
solution of \eqref{eq:min} and therefore  
${\displaystyle \lambda^1_0 \leq \frac{3}{2} \int^{1}_{-1}f(x)dx.}$  
Equality holds if and only if $u(x)=x$ is the first $S^1$ invariant 
eigenfunction. In this case, upon substitution of $u(x)=x$ into 
\eqref{eq:de} we obtain the
equivalent equation
$-f'(x)=\lambda^1_0 x$. Recalling that $f(x)$ and $f'(x)$ must satisfy 
certain boundary 
conditions forces $\lambda^1_0 = 2$ and yields the unique solution
$f(x)=1-x^2$. In other words, $g=can$.
\end{proof}

Because of Proposition 3.1 iv.), we have the immediate corollary:

\begin{cor} Let $(M,g)$ be an $S^1$ invariant Riemannian manifold of area 
$4 \pi$ which is
diffeomorphic
to $S^2$, with metric \eqref{eq:feq}. Let $ \lambda_1 $ be the first, 
non-zero, distinct
 eigenvalue for this metric, then

\begin{equation*}
\lambda_1 \leq \frac{3}{2} \int^{1}_{-1}f(x)dx
\end{equation*}
and equality holds if and only if $(M,g)$ is isometric to $(S^2,can)$. \qed
\end{cor}

\section{Spectral obstructions to isometric embeddings}

In \cite{eng2} we used the trace formula of Proposition 3.1 i.) to show that 
there exist surfaces of revolution with arbitrarily large first $S^1$ 
invariant eigenvalue. This fact, together with
 Proposition 4.1 of the last section, shows that as $\lambda_0^1$ increases 
so does 
the integral $\int^{1}_{-1}f(x)dx$. This fact
is the key to the results of this section, but first we will prove a lemma 
which gives lower bounds for our eigenvalues.

\begin{lemma}
Let $f(x)$ and $\lambda^{m}_{k}$ be defined as above then for all $m \in \N$

\begin{equation*}
\lambda^{m}_{k} > \left\{ \begin{array}{cl}
    2m \left[ \int_{-1}^{1}
 \frac{1-x^2}{f(x)} dx \right]^{-1} & 
 \mbox{if $k=0$} \\
   m|k| &  \mbox{if $k \neq 0$}
                     \end{array} \right. .
\end{equation*}

\end{lemma}

\begin{proof} From Proposition 3.1 i.) 
${\displaystyle 
\frac{1}{2} \int_{-1}^{1} \frac{1-x^{2}}{f(x)} dx = \sum^{\infty}_{j=1} 
\frac{1}{\lambda^{j}_{0}}}$ and ${\displaystyle 
\frac{1}{|k|} = \sum^{\infty}_{j=1} 
\frac{1}{\lambda^{j}_{k}}}$. 
Each of the sequences $\left\{\lambda^{j}_{k}\right\}_{j=1}^{\infty}$ is 
positive and strictly increasing so by truncating the above series after 
$m$ terms and then replacing each term with the smallest one we obtain

$$ 
\frac{1}{2} \int_{-1}^{1} \frac{1-x^{2}}{f(x)} dx > 
\frac{m}{\lambda^{m}_{0}} \hspace{.1in} \mbox{and} \hspace{.1in} 
\frac{1}{|k|} >  
\frac{m}{\lambda^{m}_{k}}.$$ This produces the desired inequalities.
\end{proof}

As was observed in \cite{eng2}, the $k=0$, $m=1$ case of this inequality, 
together with the minimal restrictions on the function $f$ is enough to 
ensure that there exist surfaces of revolution with arbitrarily large 
$\lambda^{1}_{0}$. Because of this, we can be confident that the next two
results are non-vacuous.

\begin{prop} 
Let $(M,g)$ be an $S^1$ invariant Riemannian manifold of area $4 \pi$
 which is
diffeomorphic
to $S^2$ and let $ \lambda^1_0 $ be it's first non-zero
 $S^1$
invariant eigenvalue. If $ \lambda^1_0 >3$ then $(M,g)$ cannot be 
isometrically $C^1$
embedded in $(\R^3, can)$.
\end{prop}
 
\begin{proof} By Proposition 4.1, since $ \lambda^1_0 >3$, then
 $\int^{1}_{-1}f(x)dx > 2$. Upon integrating by parts we have
 $-\int^{1}_{-1}xf'(x)dx >2$ so that 

$$2 < \left|-\int^{1}_{-1}xf'(x)dx \right| \leq \int^{1}_{-1}|x||f'(x)|dx 
\leq \max_{x \in [-1,1]}|f'(x)|.$$
So there exists $x_0 \in [-1,1]$ with $|f'(x_0)|>2$, thus, by Proposition 2.1,
precluding the possibility of an isometric embedding.
\end{proof}

Since non-embeddable metrics have some negative curvature, we have the 
immediate corollary:

\begin{cor}
Let $(M,g)$ be an $S^1$ invariant Riemannian manifold of area $4 \pi$ 
which is
diffeomorphic
to $S^2$, let $K$ be it's Gauss curvature, and let $ \lambda^1_0 $ be it's 
first
 $S^1$
invariant eigenvalue. If $ \lambda^1_0 >3$ then there exists a point $p \in
 M$
such that $K(p) < 0$. \qed
\end{cor} 

\begin{remark}
Rafe Mazzeo and Steve Zelditch have brought to our attention a recent result
of Abreu and Freitas, \cite{abr}, which is a significant improvement of 
Proposition 5.2. They prove, with the same hypothesis as Proposition 5.2 and
using the notation of this paper, 
that, for metrics isometrically embedded in $(\R^3, can)$, 
$ \lambda^j_0 < \xi_{j}^{2}/2$, for all $j$ where 
$\xi_{j}$ is a positive zero of a certain Bessel function or its 
derivative. In particular, $ \lambda^1_0 < \xi_{1}^{2}/2 \approx 2.89$. We 
have left Proposition 5.2 in the paper since its proof is so easy, and 
because the eigenvalue bound contained therein is sufficient for proving 
the main theorem (Theorem 5.5) below.
\end{remark}
 
As we allow the first $S^1$ invariant eigenvalue to increase one might 
suspect 
that, so to speak, some small eigenvalues with even multiplicity are 
``left behind". This suggests that we might find an obstruction to embedding
if the first few eigenvalues have even multiplicities. We will soon see that
even multiplicities for the first four eigenvalues will constitute such an
obstruction, but first it would be a good idea to know if metrics with this
property exist. This is the subject of:

\begin{thm} There exist metrics on $S^2$ whose first four distinct non-zero 
eigenvalues have even multiplicity.
\end{thm}

\begin{proof} To prove this theorem we will find an $S^1$ invariant 
metric of area $4 \pi$ with this property. 

By Proposition 3.1 v.), $\dim E_{\lambda_m}$ is even if and only if 
$\lambda_m$ is not an $S^1$ invariant eigenvalue, i.e. if and only if 
$ \lambda_m  \neq \lambda^j_0$ for any $j$. It is now clear that the first 
four multiplicities are even if and only if  $ \lambda_4  < \lambda^1_0$,
and, by Proposition 3.1, iv.), 
this will occur if our metric satisfies 
$ \lambda^{1}_4  
< \lambda^1_0$. Using a variational principle, as in \cite{eng3}, for 
the operator $L_4$, we obtain the upper bound:

$$\lambda^{1}_{4}\leq \frac{\int_{-1}^{1}\left[f(x)\left(\frac{du}{dx}
\right)^2+
\frac{4^2}{f(x)}u^2\right]dx}{\int_{-1}^{1}u^2dx}$$
$\forall u\in C^{\infty}
(-1,1)$ such that $u(-1)=u(1)=0$. 

Comparing this upper bound with the lower bound on $\lambda^1_0$ 
provided by Lemma 5.1., the proof of this theorem may now be reduced to 
finding a
 function $f$ and a suitable test function $u$ such that

\begin{equation}
\frac{\int_{-1}^{1}\left[f(x)\left(\frac{du}{dx}
\right)^2+
\frac{16}{f(x)}u^2\right]dx}{\int_{-1}^{1}u^2dx} < 2 \left[ \int_{-1}^{1}
 \frac{1-x^2}{f(x)} dx \right]^{-1} \label{eq:ineq}
\end{equation}

We claim that ${\displaystyle f(x)=\frac{10(1-x^2)}{1+9x^{36}}}$ and 
${\displaystyle u(x)=\sqrt{1-x^2}}$ will satisfy 
the inequality \eqref{eq:ineq}.

It is not difficult to see that ${\displaystyle 2 \left[ \int_{-1}^{1}
 \frac{1-x^2}{f(x)} dx \right]^{-1}=\frac{185}{23}>8}$ for this choice of 
$f(x)$. So the right hand side of \eqref{eq:ineq} is greater than 8.
Calculating the left hand side of \eqref{eq:ineq} for this choice of $f(x)$
and $u(x)$ yields:

\begin{eqnarray*}
 \lefteqn{\frac{\int_{-1}^{1}\left[f(x)\left(\frac{du}{dx}
\right)^2+
\frac{16}{f(x)}u^2\right]dx}{\int_{-1}^{1}u^2dx} = 
\frac{3}{4} \left[10\int_{-1}^{1}\frac{x^2}{1+9x^{36}} dx +\frac{8}{5}
\int_{-1}^{1}(1+9x^{36})dx \right]} \hspace{2.4in} \\
& & < \frac{3}{4}\left[10\cdot\frac{2}{3}+
\frac{16}{5}\cdot\frac{46}{37} \right] = \frac{1477}{185} < 8,
\end{eqnarray*} 
where the first integral in brackets has been approximated in the obvious 
way. 
Since the left hand side is less than 8 and the right hand side is greater 
than 8, the proof is finished.
\end{proof}

The proof of this theorem is hardly optimal since there are, certainly,
 many such metrics. We also believe that using a similar technique, one
should be able to find metrics whose first $m$ distinct eigenvalues have 
even multiplicity for arbitrary $m$, but we will not address these problems 
here.

\begin{thm} Let $(M,g)$ be an $S^1$ invariant Riemannian manifold which is
diffeomorphic
to $S^2$ and let $ \lambda_m $ be its $m$-th distinct eigenvalue. If 
$\dim E_{\lambda_m}$ is even for $1 \leq m \leq 4$ then  $(M,g)$ cannot be 
isometrically $C^1$
embedded in $(\R^3, can)$.
\end{thm}

\begin{proof}
Without loss of generality, we may assume the area of the metric is $4 \pi$.
As seen in the proof of Theorem 5.4, the first four eigenvalues have even 
multiplicity if and only if 
$ \lambda_4  < \lambda^1_0$. This result will then follow from Proposition 5.2 
as long as we can prove that $\lambda_4 > 3$. This is most easily 
accomplished by contradiction.

Assume $ \lambda_4 \leq 3$ so that $0<\lambda_1<\lambda_2<\lambda_3<
\lambda_4 \leq 3 $. Now each $\lambda_i$ for $1\leq i \leq 4$ must satisfy 
$\lambda_i = \lambda^l_k$ for some $k \neq 0$ and $l \geq 1$. However, by 
Lemma 5.1, $\lambda^l_k > l|k|$ so if $\lambda_i = \lambda^l_k \leq 3$ it 
must be the case that $l|k| \leq 2$. 
By Proposition 3.1 iii.) $\lambda^l_k =\lambda^l_{-k}$ so there are only  
three (possibly) distinct eigenvalues with these properties and their values
coincide with $\lambda^1_1$, $\lambda^1_2$, and 
$\lambda^2_1$.  
There are, therefore, at most three distinct values for the four 
distinct eigenvalues $ \lambda_i$ for $1 \leq i \leq 4$, but this contradicts
the pigeonhole principle.
\end{proof}

Again there is an immediate corollary:

\begin{cor}
Let $(M,g)$ be an $S^1$ invariant Riemannian manifold which is
diffeomorphic
to $S^2$, let $K$ be it's Gauss curvature, and let $ \lambda_m $ be it's 
$m$-th distinct eigenvalue. If 
$\dim E_{\lambda_m}$ is even for $1 \leq m \leq 4$ then 
there exists a point $p \in M$
such that $K(p) < 0$, \qed
\end{cor} 
and, contrapositively, a kind of partial converse for Weyl type theorems:

\begin{cor} 
A surface of revolution which is diffeomorphic to $S^2$ and 
isometrically embedded in $(\R^3, can)$ has the property that at least 
one of its 
first four non-zero distinct eigenvalues has odd multiplicity. \qed
\end{cor}

We observe that this is a property which these metrics share with those of 
constant
positive curvature.

\section{Remarks on classical surface theory}

In this final section we leave behind the question of embeddability and 
focus our attention on the way in which Corollary 5.6 can be viewed as an 
extension of one of the corollaries of the Gauss-Bonnet theorem.

Let $(M,g)$ be any compact, orientable, boundaryless surface with metric $g$.
We recall that the Euler characteristic, $\chi (M)$, and curvature $K$ are 
related by the Gauss Bonnet Theorem:

$$2\pi \chi (M)= \int_M K,$$ 
so that one has the well known result:

\begin{prop} If $\chi (M) \leq 0$ then there exists a point $p \in M$ such 
that $K(p) \leq 0$. \qed
\end{prop}

Via the Hodge-DeRham isomorphism, one can restate the Gauss-Bonnet
theorm as follows:

\vspace*{.1in}

{\em Let $\lambda_{q,j}$ be the $j$-th distinct eigenvalue of the 
Laplacian acting
on $q$-forms and $E_{\lambda_{q,j}}$ its ``eigenspace" (this vectorspace may
consist of the zero vector only). Then}

\begin{equation}
\frac{1}{2\pi}\int_M K=2-\dim E_{\lambda_{1,0}}. \label{eq:gb}
\end{equation}

\vspace*{.1in} 

Of course $\dim E_{\lambda_{1,0}}$ is simply twice the genus of the surface
 since 
$\lambda_{1,0}=0$. But this form of the Gauss-Bonnet formula does allow us
to observe that:
{\em If $\dim E_{\lambda_{1,0}}$ is even (this is automatic) and positive 
then 
there exists a point $p \in M$ such that $K(p) \leq 0$.}

In case $\dim E_{\lambda_{1,0}}>0$, $M$ is not a sphere. So these results
 tell us 
how to get some non-positive curvature by adding handles to the sphere.

If we don't want to add handles to the sphere, it is Corollary 5.6 which 
tells us, at least in the case of surfaces of revolution,
 how to obtain some negative curvature 
by changing the dimension of the euclidian space into which it embeds. 

Collecting the forgoing ideas together, one can state a result which gives
a unified, if not quite complete, answer to the question of the existence of
non-positive curvature, in other words: a generalization of Proposition 6.1 
which 
includes 
surfaces with Euler characteristic 2.

\begin{cor} Let $(M,g)$ be an orientable, compact, boundaryless 
surface with metric $g$, isometry group 
$\Im (M,g)$ and 
$j$-th distinct $q$-form eigenvalue $\lambda_{q,j}$. If, for some $q \in \{ 
| \dim \Im(M,g) -1|,1 \}$, 
$\dim E_{\lambda_{q,|1-q|\cdot j}}$ is even and positive for
 all $j$ such that $1 \leq j \leq 4$, then there exists a point $p\in M$ such 
that $K(p) \leq 0$.
\end{cor}

\begin{proof} 
If $(M,g)$ satisfies the hypothesis for $q=1$ then the statement of this 
result is simply Proposition 6.1 as can be seen from Equation \eqref{eq:gb}. 
Hence there exists a point $p\in M$ such 
that $K(p) \leq 0$.

If $(M,g)$ satisfies the hypothesis for $q= | \dim \Im(M,g) -1|$, then as is
 well known, we must consider all cases with $\dim \Im(M,g) \leq 3$. If 
$\dim \Im(M,g) = 0$ or  $2$, then, again, $q=1$ and the proof is the same
as the previous case. If $\dim \Im(M,g) = 3$ then $(M,g)=(S^2, can)$ (see 
\cite{kob}, p. 46, 47)
so that
$K>0$ and constant, and the statement of this result is simply, as we already
know, that one of the first four, 0 or 2-form, eigenvalues has odd 
multiplicity (in fact they all have odd multiplicity).
Finally, if 
$\dim \Im(M,g) = 1$ then $q=0$. If, in this case, the hypothesis holds for 
$q=0$ 
only then $M$ is, topologically, the sphere and thus the
statement of this theorem reduces to Corollary 5.6.
\end{proof}

One cannot help but ponder the possibility that one can remove the, a priori,
assumption of $S^1$ invariance since, according to legend, only surfaces of
revolution would have a lot of even multiplicities anyway. Also, if we might 
hazard an even more provacative conjecture: perhaps there is a formula 
which relates integrals of geometric invariants with multiplicities of 
non-zero eigenvalues in a similar way as Formula \eqref{eq:gb}. One could 
perhaps use heat kernel asymptotics to explore this. But this would be the 
subject of another treatise.


\begin{thebibliography}{99}

\bibitem{abr} Abreu, M., Freitas, P., {\em On the invariant spectrum of 
$S^1$-invariant metrics on $S^2$}, preprint math.SP/9909181 v2.
\bibitem{al} Alexandrov, A., {\em Intrinsic geometry of convex surfaces}, 
OGIZ Moscow, 1948.
\bibitem{be} Besse, A., {\em Manifolds All of Whose Geodesics are Closed},
 Springer-Verlag, Berlin, 1978.
\bibitem{chen} Chen, B-Y., {\em A report on submanifolds of finite type}, 
Soochow J. Math. Vol. 22, No. 2, 1996, 117-337.
\bibitem{eng1} Engman, M., {\em New Spectral Characterization Theorems for 
$S^{2}$}, Pacific J. Math. Vol. 154, No. 2, 1992, 215-229.
\bibitem{eng2} Engman, M., {\em Trace Formulae for $S^1$ invariant Green's 
Operators on $S^{2}$}, manuscripta math. 93, (1997), 357-368. 
\bibitem{eng3} Engman, M., {\em Sharp bounds for eigenvalues and 
multiplicities on surfaces of revolution}, Pacific J. Math. Vol. 186, No. 1,
 1998, 29-37.
\bibitem{gr} Greene, R., {\em Metrics and Isometric Embeddings of the
 2-sphere}, J. Diff. Geom. {\bf 5} (1971), 353-356.
\bibitem{he} Heinz, E., {\em On Weyl's embedding problem}, J. Math. 
Mech.,{\bf 11}, 421-454 (1962).
\bibitem{her} Hersch, M. J., {\em Quatre propri\'{e}t\'{e}s 
isop\'{e}rim\'{e}triques de membranes sph\'{e}riques homog\`{e}nes}, C.R.
Acad.Sc.Paris, t. 270 (1970), S\'{e}rie A, 1645-48.
\bibitem{hz} Hong, J., Zuily, C., {\em Isometric embedding of the 2-sphere
 with non negative curvature in $\R^3$.}, Math. Z. 219, 323-334 (1995).
\bibitem{hw} Hwang, A., Singer, M., {\em A Momentum Construction for Circle
-Invariant K\"{a}hler Metrics}, Preprint, math.DG/9811024.
\bibitem{kob} Kobayashi, S., {\em Transformation Groups in Differential 
Geometry}, Springer-Verlag, Heidelberg, 1995.
\bibitem{li} Li, P., Yau, S.T., {\em A new conformal invariant and its 
applications to the Willmore Conjecture and the first eigenvalue of compact
surfaces}, Invent. Math. {\bf 69} (1982), 269-291.
\bibitem{ni} Nirenberg, L., {\em The Weyl and Minkowski problems in 
differential geometry in the large}, Comm. Pure Appl. Math {\bf 6},
 337-394 (1953).
\bibitem{po} Pogorelov, A.V., {\em Extrinsic geometry of convex surfaces}, 
Trans. Math. Monogr., Vol. {\bf 6}, Providence, RI; Am. Math. Soc. 1973.
\bibitem{we} Weyl, H., {\em \"{U}ber die Bestimmung einer geschlossen 
convex ...}, Vierteljahrschrift Naturforsch. Gesell, (Z\"{u}rich), 
{\bf 61}, 40-72 (1916).

\end{thebibliography}
\end{document}